\def\filename{\jobname.tex}
\def\Stand{2000.06.19/upp}

\scrollmode 
\newlinechar=`\^^J
\def\msgline#1{\immediate\write16{^^J#1}}
\msgline{Stand: \Stand? Noch aktuell?}

\magnification=\magstep1
\vsize 23 true cm   
\hsize 15.5 true cm 
\baselineskip=13.5 pt 
\parskip=3pt 
\overfullrule=0pt

\def\TOT{}
\def\LEB{}
\long\def\TOT#1\LEB{}

\def\bR{{\bf R}} 

\def\BBF#1{\expandafter\edef\csname #1#1\endcsname{%
   {\mathord{\bf #1}}}}
\BBF m
\BBF p
\BBF C
\BBF K
\BBF N
\BBF P
\BBF Q
\BBF R
\BBF T
\BBF V
\BBF Z

\def\CAL#1{\expandafter\edef\csname #1\endcsname{%
   {\mathord{\cal #1}}}}
\CAL A
\CAL C
\CAL E
\CAL F
\CAL G
\CAL H
\CAL{IH}
\CAL J
\CAL K
\CAL L
\CAL P

\def\A{{\cal A}}
\def\B{{\cal B}}
\def\E{{\cal E}}
\def\F{{\cal F}}
\def\G{{\cal G}}
\def\H{{\cal H}}

\def\K{{\cal K}}

\newread\amssym
\openin\amssym=amssym.def
\ifeof\amssym\else\input amssym.def
\def\mm{\frak m}

\fi
\closein\amssym

\def\:{\colon}

\def\b{^{\scriptscriptstyle\bullet}}

\def\quer{\overline}

\def\longto{\longrightarrow}

\def\onto{\to\kern-.8em\to}
\def\longonto{\longto\kern-.8em\longto}

\def\epsilon{\varepsilon}
\def\phi{\varphi}
\def\rho{\varrho}
\def\theta{\vartheta}

\def\wideto{\;\to\;}

\def\qed{\hfill$\scriptscriptstyle\rfloor$}
\def\sqr#1#2{{\vbox{\hrule height.#2pt
         \hbox{\vrule width.#2pt height#1pt \kern#1pt
               \vrule width.#2pt}
         \hrule height.#2pt}}}

\def\qed{\ifmmode\hbox{\sqr35}
     \else\hfill\sqr35\fi}

\def\ibull#1\par{{\parindent10pt\item{$\bullet$}#1\par}}
\def\iibull#1\par{{\parindent10pt\itemitem{$\bullet$}#1\par}}
\def\iiitem{\par\indent\indent\hangindent3\parindent\textindent}
\def\iiibull#1\par{{\parindent10pt\iiitem{$\bullet$}#1\par}}
\def\istern#1\par{{\parindent10pt\item{$*$}#1\par}}
\def\istrich#1\par{{\parindent10pt\item{--}#1\par}}

\def\MOPnl#1{\expandafter\edef\csname #1\endcsname{%
   {\mathop{\rm #1}\nolimits}}}

\MOPnl{cld}
\MOPnl{codim}
\MOPnl{coker}
\MOPnl{Coker}
\MOPnl{Ext}
\MOPnl{Hom}
\MOPnl{Ker}
\MOPnl{id}
\MOPnl{im}
\MOPnl{ker}
\MOPnl{lin}
\MOPnl{max}
\MOPnl{odd}
\MOPnl{or}
\MOPnl{pr}
\MOPnl{pt}
\MOPnl{rg}
\MOPnl{rk}
\MOPnl{Sp}
\MOPnl{Spec}
\MOPnl{st}
\MOPnl{supp}
\MOPnl{Tor}

\font\XIIbfM=cmbx12 scaled 1200
\font\XIIbf=cmbx12

\font\XIIrm=cmr12

\font\bfit=cmbxti10
\font\sc=cmcsc10

\def\Randmarke#1{\vadjust{\vbox to 0pt
                {\vss \hbox to \hsize{\hskip\hsize
                                     \quad #1\hss}\vskip3.5pt}}}
\def\Randpfeilmarke#1{\Randmarke{$\scriptscriptstyle\Leftarrow$#1}}

\newcount\annum
\annum=\year
\ifnum\annum<2000\advance\annum by -1900\fi
\def\heute{\the\annum.\ifnum\month<10 0\fi\the\month.\ifnum
           \day<10 0\fi\the\day}

\newcount\hora
\hora=\time \divide\hora by 60
\newcount\horapostnoctem
\horapostnoctem=\hora \multiply\horapostnoctem by 60
\newcount\minuta
\minuta=\time \advance\minuta by -\horapostnoctem
\def\Hora{\ifnum\hora<10 0\fi\the\hora}
\def\Minuta{\ifnum\minuta<10 0\fi\the\minuta}
\def\Uhrzeit{\Hora:\Minuta}
\msgline{Datum: \heute, Zeit: \Uhrzeit, Min. seit 0:00:
         \the\time}

\newcount\hour      \hour\time\divide\hour60
\newcount\minute    \minute-\hour\multiply\minute60\advance\minute\time
\def\now{\ifnum\hour<10 0\fi\the\hour:\ifnum\minute<10 0\fi\the\minute}

\def\diagram{\def\normalbaselines{\baselineskip20pt\lineskip3pt
\lineskiplimit3pt}\matrix}
\def\mapright#1{\smash{\mathop{\hbox to 35pt{\rightarrowfill}}\limits^{#1}}}
\def\mapleft#1{\smash{\mathop{\hbox to 35pt{\leftarrowfill}}\limits^{#1}}}
\def\mapdown#1{\Big\downarrow\rlap{$\vcenter{\hbox{$\scriptstyle#1$}}$}}

\everydisplay{\abovedisplayskip 6pt plus 3 pt minus 2 pt
\belowdisplayskip\abovedisplayskip
\belowdisplayshortskip=3pt plus 3pt minus 1pt
}

\font\sevenrm=cmr7     \font\sixrm=cmr6
\font\seveni=cmmi7     \font\sixi=cmmi6
\font\sevensy=cmsy7    \font\sixsy=cmsy6
\font\sevenbf=cmbx7    \font\sixbf=cmbx6
\font\seventt=cmtt8 scaled 875
\font\fivei=cmmi5
\font\sevenit=cmti7
\font\sevensl=cmsl8 scaled 875

\font\fivebf=cmbx5
\def\kleinst{\def\rm{\fam0\sevenrm}
  \textfont0=\sevenrm \scriptfont0=\sixrm \scriptscriptfont0=\fiverm
  \textfont1=\seveni   \scriptfont1=\sixi  \scriptscriptfont1=\fivei
  \textfont2=\sevensy   \scriptfont2=\sixsy \scriptscriptfont2=\fivesy
  \textfont\itfam=\sevenit \def\it{\fam\itfam\sevenit}%
  \textfont\slfam=\sevensl \def\sl{\fam\slfam\sevensl}%
 \textfont\ttfam=\seventt \def\tt{\fam\ttfam\seventt}%
  \textfont\bffam=\sevenbf \scriptfont\bffam=\sixbf
   \scriptscriptfont\bffam=\fivebf \def\bf{\fam\bffam\sevenbf}%
  \normalbaselineskip=7pt
  \setbox\strutbox=\hbox{\vrule height6pt depth2pt width0pt}%
  \let\sc=\sixrm \let\big=\sevenbig \normalbaselines\rm}
\def\sevenbig#1{{\hbox{$\textfont0=\eightrm\textfont2=\eightsy
    \left#1\vbox to6.5pt{}\right.\n@space$}}}

\font\eightrm=cmr8     
\font\eighti=cmmi8     
\font\eightsy=cmsy8    
\font\eightbf=cmbx8    
\font\eighttt=cmtt8    
\font\eightit=cmti8
\font\eightsl=cmsl8
\def\kleiner{\def\rm{\fam0\eightrm}
  \textfont0=\eightrm \scriptfont0=\sixrm \scriptscriptfont0=\fiverm
  \textfont1=\eighti     \scriptfont1=\sixi \scriptscriptfont1=\fivei
  \textfont2=\eightsy     \scriptfont2=\sixsy \scriptscriptfont2=\fivesy
  \textfont\itfam=\eightit \def\it{\fam\itfam\eightit}%
  \textfont\slfam=\eightsl \def\sl{\fam\slfam\eightsl}%
 \textfont\ttfam=\eighttt \def\tt{\fam\ttfam\eighttt}%
  \textfont\bffam=\eightbf \scriptfont\bffam=\sixbf
   \scriptscriptfont\bffam=\fivebf \def\bf{\fam\bffam\eightbf}%
  \normalbaselineskip=8pt
  \setbox\strutbox=\hbox{\vrule height6pt depth2pt width0pt}%
  \let\sc=\sixrm \let\big=\eightbig \normalbaselines\rm}

\def\eightbig#1{{\hbox{$\textfont0=\ninerm\textfont2=\ninesy
    \left#1\vbox to6.5pt{}\right.\n@space$}}}
\def\ninebig#1{{\hbox{$\textfont0=\tenrm\textfont2=\tensy
    \left#1\vbox to7.25pt{}\right.\n@space$}}}

\font\ninerm=cmr9
\font\ninei=cmmi9
\font\ninesy=cmsy9
\font\ninebf=cmbx9
\font\ninett=cmtt9
\font\nineit=cmti9
\font\ninesl=cmsl9
\def\klein{\def\rm{\fam0\ninerm}
  \textfont0=\ninerm \scriptfont0=\sixrm \scriptscriptfont0=\fiverm
  \textfont1=\ninei     \scriptfont1=\sixi \scriptscriptfont1=\fivei
  \textfont2=\ninesy     \scriptfont2=\sixsy \scriptscriptfont2=\fivesy
  \textfont\itfam=\nineit \def\it{\fam\itfam\nineit}%
  \textfont\slfam=\ninesl \def\sl{\fam\slfam\ninesl}%
  \textfont\ttfam=\ninett \def\tt{\fam\ttfam\ninett}%
  \textfont\bffam=\ninebf \scriptfont\bffam=\sixbf
   \scriptscriptfont\bffam=\fivebf \def\bf{\fam\bffam\ninebf}%
  \normalbaselineskip=10pt
  \setbox\strutbox=\hbox{\vrule height7pt depth3pt width0pt}%
  \let\sc=\sevenrm \let\big=\ninebig \normalbaselines\rm}

%
%
    \newcount\annotation
     \annotation=0
      \long\def\Fussnote#1{{\baselineskip=9pt
           \setbox\strutbox=\hbox{\vrule height 7pt depth 2pt width 0pt}
         \klein\global\advance\annotation by 1
         \footnote{$^{\the\annotation)}$}{#1}}}
%
%
    \newcount\tmpcomment
     \tmpcomment=0
      \long\def\Kommentar#1{{\global\advance\tmpcomment by 1\klein
       \Randpfeilmarke{\#\the\tmpcomment}{\baselineskip=9pt
           \setbox\strutbox=\hbox{\vrule height 7pt depth 2pt width 0pt}
         \footnote{$^{\#\the\tmpcomment)}$}{\it#1}}}}
%
%
\def\makeatletter{\catcode`\@=11\relax}
\def\makeatother{\catcode`\@=12\relax}
\makeatletter

\newdimen\@tempdima
\newbox\@tempboxa

\def\@height{height}
\def\@depth{depth}
\def\@width{width}

\newdimen\fboxrule
\newdimen\fboxsep

\fboxrule = .4pt
\fboxsep = 2pt

\long\def\fbox#1{\leavevmode\setbox\@tempboxa\hbox{#1}\@tempdima\fboxrule
   \advance\@tempdima \fboxsep \advance\@tempdima \dp\@tempboxa
   \hbox{\lower \@tempdima\hbox
  {\vbox{\hrule \@height \fboxrule
          \hbox{\vrule \@width \fboxrule \hskip\fboxsep
          \vbox{\vskip\fboxsep \box\@tempboxa\vskip\fboxsep}\hskip
                 \fboxsep\vrule \@width \fboxrule}%
                 \hrule \@height \fboxrule}}}}

\makeatother


\headline{\klein\hfill Combinatorial Duality and Intersection Product
 \hfill-- \folio\ }

\ \vskip 0.5 true cm

\centerline{\XIIbfM Combinatorial Duality and
Intersection Product}
\vskip 12 pt
\centerline{\XIIrm Karl-Heinz Fieseler}
\bigskip\medskip
\centerline{\XIIbf Introduction}
\medskip
In the articles [BBFK] and [BreLu$_1$] a combinatorial
theory of intersection cohomology and perverse sheaves has 
been developed on fans. In [BBFK] we have tried to present everything on an elementary level,
using only some commutative algebra and no derived categories. 
There remained two major gaps:
 
First of all the Hard Lefschetz Theorem was only conjectured
and secondly the intersection product seemed to depend
on some non-canonical choices. Meanwhile the Hard 
Lefschetz theorem has been proved 
in [Ka]. The proof relies heavily on the intersection product,
since what finally has to be shown are the Hodge-Riemann relations. 
In fact here again choices enter: The intersection product 
is induced from the intersection product on some
simplicial subdivision via a direct embedding of the 
corresponding intersection cohomology sheaves, a fact, which makes  
the argumentation quite involved. In their paper [BreLu$_2$]
Bressler and Lunts show by a detailed analysis that 
eventually all possible choices do not
affect the definition of the pairing.

Our goal here is the same, but we shall try to follow
the spirit of [BBFK] avoiding the formalism of derived categories. 
For perverse sheaves we define their dual sheaf and 
check that the intersection cohomology sheaf is self-
dual. Since, as a consequence of 
Hard Lefschetz, an isomorphism between 
intersection cohomology sheaves is unique
up to multiplication with a real number, we obtain,
after the choice of a volume form on $V$, a natural intersection 
product. In fact
that fits completely into the proof of HL according to [Ka],
which uses induction on the dimension of the fans involved. 
The HL in dimensions $ < n$ already gives us the naturality
of the intersection product in dimension $n$, in particular, it
agrees with the restriction (with respect to a direct embedding) 
of the intersection product
on any simplicial refinement of the given fan. 
 
This article is a preliminary version of a paper
in collaboration with G.Barthel, J.-P. Brasselet
and L.Kaup. The author is grateful for an invitation
to the IML at Marseille-Luminy, where this article has 
been written. 

\bigskip\medskip\goodbreak
\centerline{\XIIbf 0. Preliminaries}
\medskip\nobreak
Let $V$ be a real vector space of dimension
$n:=\dim V$ and denote $A\b:=S\b (V^*)$ the symmetric algebra
on $V^*$, with other words the algebra of real valued
polynomials on $V$, graded "topologically", i.e., 
such that $V^*=A^2$. Let $\mm:=A^{>0}$ and, for
a graded $A\b$-module $M\b$ denote
$$
\overline M\b:=A\b/\mm \otimes_{A\b} M\b
$$
its reduction mod $\mm$, a graded real vector space.
For a strictly convex
polyhedral cone $\sigma \subset V$ denote
$V_\sigma \subset V$ its linear span, and
$$
A\b_\sigma:=S\b(V_\sigma^*)
$$
again with the convention $A^2_\sigma=V_\sigma^*$.

We consider a fan $\Delta$ in $V$
as a topological space with the subfans as open subsets.
Sheaf theory
on such a space is particularly simple since the
``affine'' open fans $\langle \sigma \rangle
\le \Delta$ consisting of a cone $\sigma$ and its 
proper faces
form a basis of the fan topology
whose elements can not be
covered by strictly smaller open sets. In fact,
let $(F_{\sigma})_{\sigma \in \Delta}$
be a collection of abelian groups, say, together
with ``restriction'' homomorphisms
$\rho^\sigma_{\tau} \: F_{\sigma} \to F_{\tau}$
for $\tau \le \sigma$, i.e., we require
$\rho^\sigma_\sigma = \id$ and
$\rho^\tau_\gamma \circ \rho^\sigma _\tau =
\rho^\sigma_\gamma$
for $\gamma \le \tau \le \sigma$.
Then there is a
unique sheaf~$\F$ on
the fan space~$\Delta$ such that its group of
sections $\F(\sigma)  :=
\F\bigl(\left< \sigma \right> \bigl)$ agrees
with $F_{\sigma}$. The sheaf $\F$ is flabby if
and only if each restriction map $\rho^\sigma_{\partial \sigma}
\: \F(\sigma) \to \F(\partial \sigma)$ is
surjective.  \par

In particular, we consider the sheaf~$\A\b$
of graded polynomial algebras on~$\Delta$
determined by $\A\b(\sigma)  :=  A\b_\sigma$,
the homomorphism
$\rho^\sigma_{\tau} \: A\b_\sigma \to
A\b_\tau$ being the restriction of functions
on $\sigma$ to the face $\tau \le \sigma$.
The set of sections
$\A\b(\Lambda)$ on a subfan $\Lambda \le
\Delta$ constitutes the algebra of
\hbox{($\Lambda$-)} piecewise polynomial
functions on $|\Lambda|$ in a natural way.
\par
For notational convenience, we often write
$$
  F\b_\Lambda  :=  \F\b(\Lambda)
  \qquad\hbox{and}\qquad
  F\b_\sigma  :=
  \F\b(\sigma)\; ;
$$
more generally, for a pair of subfans
$(\Lambda, \Lambda_0)$ with $\Lambda_0 \le \Lambda$, we define
$$
  F\b_{(\Lambda, \Lambda_0)}  :=
  \ker (\rho^\Lambda_{\Lambda_{0}} \:
  F\b_\Lambda \longto F\b_{\Lambda_0})\,,
$$
the submodule of sections on $\Lambda$
vanishing on $\Lambda_0$.  In particular, for
a purely $n$-dimensional fan~$\Delta$, we
obtain in that way the module
$$
  F\b_{(\Delta, \partial \Delta)}  :=
  \ker (\rho^\Delta_{\partial \Delta} \:
  F\b_\Delta \longto F\b_{\partial \Delta})
$$
of sections over $\Delta$ with ``compact
supports''. Here the boundary fan 
$\partial \Delta$ is the union of the
affine fans $\langle \tau \rangle$
with the cones $\tau$ of dimension $n-1$ being
a facet of precisely one $n$-dimensional cone
$\sigma \in \Delta$. Furthermore we write
simply $\partial \sigma$ instead of
$\partial \langle \sigma \rangle=
\langle \sigma \rangle \setminus \{ \sigma \}$.

Finally let $f:V \longrightarrow W$ be a linear map
inducing a map of fans between the fan 
$\Lambda$
in $V$ and the fan $\Delta$ in $W$, i.e.
for every cone $\sigma \in \Lambda$ there is a cone $\tau
\in \Delta$ with $f(\sigma) \subset \tau$.
Denote $\A\b$
resp. $\B\b$ the corresponding sheaves of piecewise linear
polynomials, let $\F\b$ resp. $\G\b$ be a sheaf of
$\A\b$- resp. $\B\b$-modules. Then the direct image
$f_*(\F\b)$ is the $\B\b$-module sheaf defined by
$$
f_*(\F\b) (\tau):=F\b_{f^{-1}(\tau)}
$$
with the subfan
$f^{-1}(\tau):=\{ \sigma \in \Lambda; f(\sigma)
\subset \tau \}$. We mention here in particular
the case, where $f=\id_V: V \longrightarrow V$
and $\Lambda:= \hat \Delta$ is a refinement of
$\Delta$.
And the pull back
$f^*(\G\b)$ is given by
$$
f^*(\G\b)(\sigma):=A\b_\sigma \otimes_{B\b_\tau}
G\b_\tau\ ,
$$ 
where $\tau \in \Lambda$ is the smallest cone
containing $f(\sigma)$.

\bigskip\medskip\goodbreak
\centerline{\XIIbf 1. Perverse Sheaves on a Fan}
\medskip
\noindent
{\bf 1.1 Definition:} {\it A
{\bfit perverse\/} sheaf on a
fan $\Delta$ is a {\it flabby\/}
sheaf~$\F\b$ of graded $\A\b$-modules such
that, for each cone $\sigma \in \Delta$, the
$A\b_{\sigma}$-module $F\b_{\sigma}$ is {\it
finitely generated and free}.}

\smallskip
Perverse sheaves are built up from simple objects,
which are indexed by the cones $\sigma \in \Delta$:

\medskip
\noindent
{\bf 1.2 Simple Perverse Sheaves:\/ }
For each cone $\sigma \in \Delta$, we
construct inductively a ``simple''  sheaf
${}_{\sigma}\L\b$ on~$\Delta$ as follows:
For a cone  $\tau \in \Delta$ with
$\dim \tau \le \dim \sigma$, we
set
$$
  {}_{\sigma}L\b_{\tau} \; := \;
  {}_{\sigma}\L\b(\tau) \; := \;
  \cases{ A\b_{\sigma} & if $\tau = \sigma$, \cr
             0        & otherwise.\cr }
$$
If ${}_{\sigma}\L\b$ has been defined
on $\Delta^{\le m}$ for some $m \ge \dim
\sigma$, then for each $\tau \in
\Delta^{m+1}$, we choose a linear map
${}_{\sigma}\quer L\b_{\partial \tau}
\buildrel s \over
\longrightarrow {}_{\sigma}L\b_{\partial \tau}$,
a section of the reduction map
${}_{\sigma}L\b_{\partial \tau}
\longrightarrow {}_{\sigma}\quer L\b_{\partial \tau}$,
and set 
$$
  {}_{\sigma}L\b_{\tau} \; := \;
  A\b_{\tau} \otimes_{\RR} \quer
  {}_{\sigma}\quer L\b_{\partial \tau}
  \quad
$$
and define the restriction map 
$ \rho^\tau_{\partial\tau}$ as the map obtained from
$s$ by extension of coefficients.
\smallskip

Let us collect some useful facts about
these sheaves, proved in [BBFK] and [BreLu$_1$]:

\medskip\noindent
{\bf 1.3 Remark:} i) The sheaf $\F\b :=
{}_{\sigma}\L\b$ is perverse; it is determined by the
following properties:
\itemitem{a)} $\quer F\b_{\sigma}
\cong \RR\b$,
\itemitem{b)} for each cone $\tau
\ne \sigma$, the reduced restriction map
$\quer F\b_{\tau} \to \quer F\b_{\partial
\tau}$ is an isomorphism.  \par\smallskip

\noindent
ii) For the zero cone~$o$, the simple sheaf
$$
\E\b:={}_{o}\L\b
$$ 
is called the intersection cohomology sheaf or
minimal extension sheaf
of~$\Delta$. For a quasi-convex fan $\Delta$, i.e.
such that $\Delta$ is purely $n$-dimensional
and $|\partial \Delta|$ a real homology manifold, we may then 
define its intersection cohomology
$$
IH\b (\Delta):= \overline E\b_\Delta\ .
$$  

\noindent
iii) The sheaf ${}_{\sigma}\L\b$ vanishes outside
the star
$$
\st_{\Delta}(\sigma):=\bigcup_{\tau \ge \sigma}
\langle \tau \rangle\ .
$$
Denote $f:V \longrightarrow W:=V/V_\sigma$ the
projection. It
induces a map of fans $
\st_{\Delta}(\sigma) \longrightarrow
\Delta_\sigma:=f(\st_{\Delta}(\sigma))$.
Then
$$
{}_{\sigma}\L\b \cong
f^*({}_{\Delta_\sigma}\E\b)
$$
with the intersection cohomology sheaf 
${}_{\Delta_\sigma}\E\b$ of the "transversal" 
fan $\Delta_\sigma$.

\smallskip\noindent
The following decomposition theorem has been proved in
[BBFK]:

\medskip\noindent
{\bf 1.4 Decomposition Theorem:}
{\it Every perverse sheaf $\F\b$ on~$\Delta$
admits a direct sum decomposition
$$
  \F\b \;\cong\;\;
  \bigoplus_{\sigma \in \Delta}
  \bigl({}_{\sigma}\L\b \otimes_{\RR}
  K\b_\sigma\bigr) 
$$
with $K\b_\sigma  :=  K\b_\sigma(\F\b)  :=
\ker\,(\,\quer \rho^{\sigma}_{\partial\sigma}
\: \quer F\b_{\sigma} \to
\quer F\b_{\partial\sigma})$, a finite
dimensional graded vector space.}
\smallskip
\noindent
{\bf 1.5 Example:} Let $\hat \Delta$ be a refinement
of the fan $\Delta$ and consider the identity map
$\pi:=\id_V$ as a map between the fans $\hat \Delta$
and $\Delta$ with intersection cohomology sheaves
$\hat \E\b$ and $\E\b$. Then 
$\pi_*(\hat \E\b)$ is a perverse sheaf and
we obtain a decomposition
$$
  \pi_*({\hat \E}\b) \;\cong\;\;
  \E\b \oplus \bigoplus_{\sigma \in
  \Delta^{> 1}} \kern -3pt
  {}_{\sigma}\L\b \otimes_{\RR} K\b_\sigma
$$
of $\A\b$-modules, where the~$K\b_\sigma$ are
(positively) graded vector spaces, and the ``correction terms''
are supported on the closed subset $\Delta^{> 1}:=
\Delta \setminus \Delta^{\le 1}$, where
$\Delta^{\le 1} \subset \Delta$ denotes the 1-skeleton
of $\Delta$.
\bigskip
\bigskip\medskip\goodbreak
\centerline{\XIIbf 2. The Dual of a Perverse
Sheaf}
\medskip\nobreak\noindent
{\bf 2.1 Definition:} {\it Let $\F\b$ be a perverse sheaf.
The dual sheaf $D\F$ is defined by
$$
(D\F)\b_\sigma:= \Hom (F\b_{(\sigma, \partial \sigma)},
A\b_\sigma) \otimes \det V_\sigma^* . 
$$
For a face $\tau < \sigma$ the restriction homomorphism 
$(D\F)\b_\sigma \longrightarrow (D\F)\b_\tau$ is first defined
in the case of a facet $\tau <_1 \sigma$ and in general obtained as the
composition of such maps depending a priori on the choice
of an ascending chain $\tau_0:= \tau <_1 \tau_1 <_1 
...<_1 \tau_r:=\sigma$.} 
\bigskip
\noindent
{\bf The definition of the restriction map 
$\rho^{\sigma}_\tau$ for a facet
$\tau$ of the cone $\sigma$}:
Fix a linear form $h \in V_\sigma^*$ having kernel
$\ker (h)=V_\tau$. We define maps 
$$
\phi_h: \Hom ( F\b_{(\sigma, \partial \sigma)}, 
A\b_\sigma)
\longrightarrow 
\Hom ( F\b_{(\tau, \partial \tau)}, 
A\b_\tau)\ ,\ 
\psi_h: \det V_\sigma^* \longrightarrow
\det V_\tau^*
$$
of degree 2 resp. -2,
satisfying
$$
\phi_{\lambda h}=\lambda \phi_h\ , \ 
\psi_{\lambda h}=\lambda^{-1} \psi_h\ .
$$
for any non-zero real number $\lambda \in \bR$.

\noindent
Thus there is a well defined map
$$
r^\sigma_\tau := \phi_h \otimes \psi_h:
$$
$$
(D\F)\b_\sigma= \Hom ( F\b_{(\sigma, \partial \sigma)}, 
A\b_\sigma) \otimes \det V_\sigma^*
\longrightarrow
(D\F)\b_\tau= \Hom ( F\b_{(\tau, \partial \tau)}, 
A\b_\tau) \otimes \det V_\tau^*,
$$
homogeneous of degree 0.
The map $\psi_h$ associates to $h \wedge \eta
\in \det V_\sigma^*$ the image of
$\eta \in \bigwedge^{s-1}V_\sigma^*$ (with
$s:= \dim \sigma$) under the natural map
$
\bigwedge^{s-1}V_\sigma^*
\longrightarrow
\bigwedge^{s-1}V_\tau^*=\det V_\tau^*. 
$

In order to
define~$\phi_h$, we
use three exact sequences, starting with
$$
0 \wideto
F\b_{(\sigma,\partial \sigma)} \wideto
F\b_\sigma \wideto  F\b_{\partial \sigma}
\wideto 0 \;.
\leqno(2.1.1)
$$
The second one is composed of the
multiplication with $h$ and the
projection onto the cokernel:
$$
0 \wideto A\b_\sigma~
{\buildrel \mu_h \over \longto}~
 A\b_\sigma  \wideto
 A\b_\tau \wideto 0\;.
\leqno(2.1.2)
$$
Eventually the subfan $\partial_\tau \sigma
:=  \partial\sigma \setminus \{\tau\}$ of
$\partial\sigma$ yields the exact sequence
$$
0 \wideto F\b_{(\tau, \partial \tau)}
\wideto F\b_{\partial \sigma}
\wideto F\b_{\partial_\tau  \sigma}
\wideto 0\;.\leqno(2.1.3)
$$
The associated $\Hom$-sequences provide
a diagram
$$
\klein
\def\longto{\kern-6pt\longrightarrow\kern-6pt}
\diagram{
    &   &   &   & \Ext(F_{\partial_\tau\sigma}\b,
A\b_\sigma)
                    &   &   \cr 
{\kern-50pt(2.1.4)}
    &   &   &   & \mapdown{}
                    &   &   \cr 
  \Hom(F\b_\sigma, A\b_\sigma) & \longto &
          \Hom(F\b_{(\sigma,\partial\sigma)}, 
A\b_\sigma) &
              \buildrel \alpha\over\longto &
                  \Ext(F\b_{\partial\sigma}, 
A\b_\sigma)
                    &   &   \cr 
    &   &   &   & \mapdown{\beta}
                    &   &   \cr 
  \Hom( F_{(\tau,\partial\tau)}\b,  
A\b_\sigma ) & \longto &
          \Hom( F_{(\tau,\partial\tau)}\b,  
A_\tau\b) &
              \buildrel \gamma\over\longto &
                  \Ext(F_{(\tau,\partial\tau)}\b,  
A\b_\sigma)
                    & \longto &
\Ext(F_{(\tau,\partial\tau)}\b, A\b_\sigma )
\cr 
}
$$
with $\Ext=\Ext^1_{A\b}$. We show that~$\gamma$
is an isomorphism; we then may set
$$
  \phi_h :=
  \gamma^{-1}\circ\beta\circ\alpha \,.
$$
Indeed the rightmost arrow in the bottom row
is the zero homomorphism, since it is induced by
multiplication with $h$, which annihilates
$F\b_{(\tau, \partial \tau)}$. On the other hand,
the $A\b_\tau$-module $F\b_{(\tau, \partial \tau)}$
is a torsion module over $A\b_\sigma$, so that
$\Hom(F_{(\tau,\partial\tau)}\b, A\b_\sigma
)$
vanishes.
\par

An explicit description of $\phi_h$
is as follows: For a homomorphism $f \:
F\b_{(\sigma, \partial \sigma)} \to
A\b_\sigma$, the map
$\phi_h (f) : F\b_{(\tau, \partial \tau)} \to
 A\b_\tau$ acts in the following way: Extend
a section
$g \in F\b_{(\tau, \partial \tau)}$ first trivially
to the whole boundary $\partial \sigma$ and then
to a
section $\hat g \in F\b_\sigma$; the 
map $\phi_h (f)$ now sends
$g$ to
$f(h \hat g)|_\tau$. In particular we see
that, if $\tau, \tau' <_1 \sigma$ are facets
of $\sigma$
meeting in a common facet $\gamma <_2 \sigma$,
then
$$
r^{\tau}_\gamma \circ r^{\sigma}_\tau
+r^{\tau'}_\gamma \circ r^{\sigma}_{\tau'}=0\ .
$$
Now fix for any cone $\in \Delta$ an orientation,
such that on $n$-dimensional cones the resulting map 
is constant. For a facet $\tau <_1 \sigma$
of a cone $\sigma \in \Delta$ denote
$\epsilon^{\sigma}_\tau= \pm 1$ the corresponding
transition coefficient. Finally the restriction map
$$
\rho^{\sigma}_\tau:
(D\F)\b_\sigma \longrightarrow 
(D\F)\b_\tau
$$
is defined as
$$
\rho^{\sigma}_\tau:=\epsilon^{\sigma}_\tau 
\cdot
r^{\sigma}_\tau\ .
$$ 
Then we have in the above situation
$$
\rho^{\tau}_\gamma \circ \rho^{\sigma}_\tau
=\rho^{\tau'}_\gamma \circ \rho^{\sigma}_{\tau'}\ .
$$
Hence in order to construct 
$\rho^{\sigma}_\tau$ for an arbitrary face
$\tau < \sigma$ we may choose 
an ascending chain $\tau_0:=\tau <_1 \tau_1 <_1...
<_1 \tau_r:=\sigma$ and define 
$\rho^\sigma_\tau$ as the composition of the 
restriction maps $\rho^{\tau_{i+1}}_{\tau_i}, i=0,...
, r-1$.
\bigskip
\noindent
{\bf 2.2 Theorem:} {\it The dual sheaf $D\F$ of a perverse
sheaf $\F\b$ is again perverse.}
\bigskip
\noindent
{\bf Proof}: Since according to Cor. 4.12 in [BBFK]
the $A\b_\sigma$-module $F\b_{(\sigma, \partial \sigma)}$
is free,
all we have to prove is that for every cone
$\sigma \in \Delta$ the restriction 
map
$$
(D\F)\b_\sigma \longrightarrow (D\F)\b_{\partial \sigma}
$$
is onto.
In order to do that we have to understand 
$(D\F)\b_{\partial \sigma}$. We may assume
$\dim \sigma = n$ and fix a line
$\ell \subset V$ intersecting $\buildrel \circ \over 
\sigma$ and denote $\pi: V \longrightarrow
W:=V/ \ell$ the quotient projection, $\Lambda:=
\pi (\partial \sigma)$. Furthermore
let $B\b:=S\b(W^*)$ and
$\G\b:= \pi_* (\F\b|_{\partial \sigma})$. 
Then there is a natural
isomorphism 
$$
(D\G)\b_\Lambda \cong (D\F)\b_{\partial \sigma}
$$
of $B\b$-modules, while
$$
(D\G)\b_\Lambda \cong \Hom_{B\b} (G\b_\Lambda, B\b) \otimes 
\det W^*\ .
$$   
Since that isomorphism plays an essential role
in the definition of the intersection product
we state it in a theorem:
\bigskip
\noindent
{\bf 2.3 Theorem:} {\it For a perverse sheaf $\F\b$ 
on a quasi-convex fan $\Delta$
there is a natural isomorphism}
$$
(D\F)\b_\Delta \cong \Hom_{A\b} (F\b_{(\Delta,
\partial \Delta)},
 A\b) \otimes 
\det V^*\ .
$$   
\bigskip
\noindent
Let us assume 2.3 for the moment. Denote $\hat \sigma$
the fan obtained from $\langle \sigma \rangle$ by adding the new
ray $\rho:=\ell \cap \sigma$. The projection
$\pi: V \longrightarrow W$ induces a map of fans
$\pi: \hat \sigma \longrightarrow \Lambda$, denote
$\H\b:=\pi^*(\G\b)$. We have an isomorphism
$$
H\b_{\hat \sigma}
\cong A\b \otimes_{B\b} G\b_\Lambda
\cong A\b \otimes_{B\b} F\b_{\partial \sigma}\ ,
$$
as well as the Thom isomorphism
$$
H\b_{(\hat \sigma, \partial \hat \sigma)}
= g H\b_{\hat \sigma}
\cong
H\b_{\hat \sigma}[-2]
$$
with a non-trivial function $g \in
A^2_{(\hat \sigma, \partial \hat \sigma)}$ - in fact
$g$ is unique up to a non-zero scalar 
multiple. The function $g$ induces
an isomorphism
$\det V^* \cong \det V^*_\tau,
g_\tau \wedge \eta
\mapsto \eta|_{V_\tau}$ with
$g_\tau:=g|_{\hat \tau} \in V^*, \hat \tau:=
\tau + \rho$, such that
the isomorphism
$\det V^* \cong \det V^*_\tau
\cong \det W^*$ does not depend on the facet
$\tau<_1 \sigma$.
So in order to get rid of the factors 
$\det (..)$ it is sufficient to
choose an isomorphism $\bR \cong \det V^*$.

\noindent
Now, given a section
$\in (D\F)\b_{\partial \sigma}$, or with other words,
a $B\b$-module homomorphism
$\beta:F\b_{\partial \sigma} \longrightarrow B\b$ 
an inverse image $\alpha \in (D\F)\b_{\sigma}$ with 
respect to the restriction map $\rho^\sigma_
{\partial \sigma}$ can be found as 
follows: The
map $\alpha: F\b_{(\sigma, \partial \sigma)}
\longrightarrow A\b$
is the composition 
$$
F\b_{(\sigma, \partial \sigma)}
\buildrel i \over  \longrightarrow
H\b_{(\hat \sigma, \partial \hat \sigma)}
\buildrel g^{-1} \cdot ... \over
\longrightarrow 
H\b_{\hat \sigma} \cong A\b \otimes_{B\b} 
F\b_{\partial
\sigma}
\buildrel A\b \otimes \beta \over
\longrightarrow A\b\ .
$$ 
Here the first map $i$ is obtained as the restriction
of a map $j:F\b_{\sigma}
\longrightarrow H\b_{\hat \sigma}$
factoring the restriction $F\b_{\sigma}
\longrightarrow F\b_{\partial \sigma}$
as
$$
F\b_{\sigma}
\buildrel j \over 
\longrightarrow 
H\b_{\hat \sigma}
\longrightarrow
H\b_{\partial \hat \sigma} \cong
F\b_{\partial \sigma}
$$
Such a map clearly exists, since $F\b_\sigma$ is a free
$A\b$-module. \qed
\bigskip
\noindent
{\bf Proof of 2.3:} Denote $h= h_1 \cdot...\cdot h_r \in A\b$ the product 
of all the linear forms having a $V_\tau, \tau \in
\Delta^{n-1} \setminus \partial \Delta$ as its
kernel, where different cones $\tau$ may give rise
to the same subspace $V_\tau$, but the corresponding factor
in $h$ should be simple. After having fixed 
an isomorphism
$\bR \cong \det V^*
\buildrel
\psi_{h_i}\ \cong
\over
\longrightarrow
\det V_i^*
$
with $V_i:= \ker (h_i)$ we may, for
a cone $\sigma \in \Delta^{\ge n-1}$,
replace $(D\F)\b_\sigma$
with $\Hom (F\b_{(\sigma, \partial \sigma)}, 
A\b)$ and restriction maps with 
$\pm \phi_{h_i}$.

For an $n$-dimensional
cone $\sigma \in \Delta$ the trivial extension map
$F\b_{(\sigma, \partial \sigma)} \longrightarrow 
F\b_\Delta$ induces a map
$$
\Hom (F\b_\Delta, A\b)
\longrightarrow
(D\F)\b_\sigma=\Hom (F\b_{(\sigma, \partial \sigma)}, 
A\b)
$$
factorizing (uniquely) through
$(D\F)\b_\Delta$: Check that for a cone $\tau \in
\Delta^{n-1}$ the map
$$
\Hom (F\b_\Delta, A\b)
\longrightarrow
(D\F)\b_\sigma
\longrightarrow
(D\F)\b_\tau
$$
does not depend on the choice of the cone $\sigma
\in \Delta^n, \sigma > \tau$. The resulting map
$$
\Hom (F\b_\Delta, A\b)
\longrightarrow
(D\F)\b_\Delta
$$
is injective, since the submodule
$$
F\b_{(\Delta, \Delta^{n-1})}
\cong
\bigoplus_{\sigma \in \Delta^n} F\b_{(\sigma,
\partial \sigma)}
\subset F\b_\Delta
$$
has maximal rank. On the other hand any section
$\in (D\F)\b_\Delta$ 
given by a collection of homomorphisms
$\psi_\sigma:F\b_{(\sigma, \partial \sigma)}
\longrightarrow A\b, \sigma \in \Delta^n$,
induces a homomorphism
$\psi := \sum_{\sigma \in \Delta^n}
\psi_\sigma: F\b_{(\Delta, \Delta^{n-1})} \longrightarrow
A\b$. Now $\psi$ extends uniquely
to a homomorphism $\hat \psi:
F\b_{(\Delta, \partial \Delta)}
\longrightarrow h^{-1}A\b$ and we show 
that already $\hat \psi(
F\b_{(\Delta, \partial \Delta)})
\subset A\b$ resp. 
$\psi(h
F\b_{(\Delta, \partial \Delta)})
\subset hA\b$. 
Let $f \in F\b_{(\Delta, \partial \Delta)}$
and $f_\sigma:=f|_\sigma$. Then 
$$
\psi (hf)= \sum_{\sigma \in \Delta^n}
\psi_\sigma (hf_\sigma) \in A\b
$$
is divisible by each linear factor $h_i \in A^2$.
Write $h=g_i h_i$. If
$\sigma \cap V_i$ is not a facet of
$\sigma$, we have
$g_i f_\sigma \in F\b_{(\sigma, \partial \sigma)}$
and hence $\psi_\sigma (hf_\sigma)=h_i \psi_\sigma 
(g_if_\sigma)
\in hA\b$. The remaining cones can be grouped
in pairs of different cones $\sigma, \sigma'$ such that
$\tau:=V_i \cap \sigma=V_i \cap \sigma'$ is a facet of both
$\sigma$ and $\sigma'$, and we are done if we can show
that $h_i$ divides $\psi_\sigma (hf)+
\psi_{\sigma'} (hf)$. But
$$
\psi_\sigma (hf)|_{V_i}= \phi_{h_i} (\psi_\sigma)
((g_if)|_{\tau})\ ,
$$
and the analogous formula holds for $\sigma'$ 
instead of $\sigma$. On the other hand we know that
$$
\phi_{h_i} (\psi_\sigma)=-\phi_{h_i} 
(\psi_{\sigma'})\ ,
$$
since the $\psi_\sigma, \sigma \in \Delta^n$, patch 
together to a section of $D\F$. That gives the desired
divisibility. (Remember here, that the transition coefficients
$\epsilon^\sigma_\tau$ enter in the definition of
the restriction map $\rho^\sigma_\tau$!)
\qed
\bigskip
\noindent
In order to see that the dual sheaf $D\F$ of a simple
perverse sheaf again is simple, we need biduality. In fact, it 
is an immediate consequence of the following
\bigskip
\noindent
{\bf 2.4 Proposition}: {\it For a cone $\sigma \in
\Delta$ and a perverse sheaf $\F\b$ we have
$$
(D\F)\b_{(\sigma, \partial \sigma)}=
\{ \psi|_{F\b_{(\sigma, \partial \sigma)}};\psi \in  
\Hom (F\b_{\sigma}, A\b_\sigma) \} \otimes 
\det V^*_\sigma  
$$}
\bigskip
\noindent
{\bf Proof}: The inclusion "$\supset$" is obvious
from the definition of the restriction maps for the 
dual sheaf $D\F$.

\noindent "$\subset$":
Denote $\tau_1,...,\tau_r <_1 \sigma$ 
the facets of the cone $\sigma$, let $V_{\tau_i}
= \ker (h_i), h:=h_1 \cdot... \cdot h_r$
and use the same identifications as in the proof of 
2.3. It suffices
to show that for $\psi \in (D\F)\b_{(\sigma, \partial \sigma)}$
we have $\psi (hF\b_\sigma) \subset
h A\b$. Let $\psi \in \Hom 
(F\b_{(\sigma, \partial \sigma)}, A\b_\sigma)
\cap (D\F)\b_{(\sigma, \partial \sigma)}$. Take a section
$f \in F\b_\sigma$. Write $h=g_ih_i$. Then
$g_i f$ is an extension of $(g_if)|_{\tau_i}$
as required in the definition of $\psi|_{\tau_i}=0$.
Hence $\psi (h_i(g_i f)) \in h_i A\b$ for all
$i=1,...,r$ resp., since the
$h_i,i=1,...,r$ are pairwise relatively prime,
$\psi (hf) \in hA\b$. \qed
\bigskip
\noindent
{\bf 2.5 Theorem}: {\it Let $\F\b$ be a perverse sheaf on the fan 
$\Delta$.
There is a natural isomorphism
$$
\F\b \buildrel \cong \over \longrightarrow D(D\F)\ .
$$
of perverse sheaves.} 
\bigskip
\noindent
{\bf Proof}: After the choice of an isomorphism
$\det V^* \cong \bR$ the inclusion
$$
(D\F)\b_{(\sigma, \partial \sigma)}
\subset 
(D\F)\b_{\sigma}
$$
can be interpreted as the (injective) restriction map
$$
\Hom (F\b_\sigma, A\b) \buildrel
res \over \hookrightarrow
\Hom (F\b_{(\sigma, \partial \sigma)}, A\b)\ .
$$
Hence $D(D\F)\b_\sigma \cong 
(F\b_\sigma)^{**}$
and
our sheaf isomorphism over $\sigma
\in \Delta$ can be defined as the biduality isomorphism
$$
\beta_\sigma :
 F\b_\sigma \longrightarrow (F\b_\sigma)^{**}
$$
for the free $A\b$-module $F\b_\sigma$. We leave it to the reader to
check that the above family $(\beta_\sigma)_{\sigma
\in \Delta}$ of isomorphisms defines a 
sheaf
homomorphism.
\qed
\bigskip
\noindent
{\bf 2.6 Corollary:} {\it For the simple perverse
sheaf ${}_{\sigma}\L\b$ with a cone $\sigma \in \Delta$
we have 
$$
D ({}_{\sigma}\L) \cong {}_{\sigma}\L\b
\otimes \det V_\sigma^*\ ,
$$
in particular, the intersection cohomology sheaf
$\E\b:={}_{o}\L\b$ is self dual:
$$
D\E \cong \E\b\ .
$$}
\bigskip
\noindent
{\bf Proof}: We have obviously
supp$(D\F)=$supp$(\F\b)$. Furthermore
$D(\F\b \oplus \G\b) \cong D\F \oplus
D\G$. Now use biduality in order to see that
simple sheaves are mapped by the duality functor to 
simple ones.
Finally check that 
$D ({}_{\sigma}\L)\b_\sigma \cong {}_{\sigma} L\b_\sigma
\otimes \det V_\sigma^*$.
\qed
\bigskip
\bigskip\medskip\goodbreak
\centerline{\XIIbf 3. The Intersection Product}
\medskip\nobreak\noindent
Fix a volume form
$\omega \in \det (V^*)$
resp. an isomorphism $\det (V^*) \cong \bR$. For
a quasi-convex fan $\Delta$ it induces
an isomorphism
$$
(D\E)\b_\Delta \cong \Hom (E\b_{(\Delta, \partial \Delta)},
A\b) \otimes \det V^* 
\cong 
\Hom (E\b_{(\Delta, \partial \Delta)},
A\b[-2n])\ .\leqno (D(\omega)) 
$$
According to Cor.2.6 there is an isomorphism
$$
\E\b \cong D\E \leqno (DI)
$$
thus giving rise to a map
$$
E\b_\Delta \buildrel \cong \over
\longrightarrow 
\Hom (E\b_{(\Delta, \partial \Delta)},
A\b[-2n])
$$
resp. a dual pairing
$$
E\b_\Delta \times E\b_{(\Delta, \partial \Delta)}
\longrightarrow A\b[-2n]\ .\leqno (DP)
$$
The naturality of this "intersection product" is a 
consequence of the following proposition 3.1. In fact,
this is
the only place, where the
Hard Lefschetz theorem, cf. [Ka], enters.
\bigskip
\noindent
{\bf 3.1 Theorem}: {\it Every homomorphism $\E\b
\to \F\b$ of degree 0 between two copies $\E\b,\F\b$ of the 
intersection cohomology sheaf
on $\Delta$ is determined by the homomorphism
$\RR\b \cong E\b_o \to F\b_o \cong \RR\b$.}
\medskip
\noindent
For the proof we need the following vanishing lemma:
\medskip
\noindent
{\bf 3.2 Lemma:} {\it For the intersection cohomology
sheaf $\E\b$ on a fan $\Delta$ and a non-zero cone $\sigma \in
\Delta$ we have
$$
\overline E_\sigma^q = 0 \ , 
\ \ q \ge \dim \sigma \ \ \ \ \ {\it as\ well\ as}
\ \ \ \ \ 
\overline E_{(\sigma, \partial \sigma)}^q = 0 \ ,
\ \ 
q \le \dim \sigma\ . 
$$
In particular even
$$
E_{(\sigma, \partial \sigma)}^q = 0\ , \ \ 
q \le \dim \sigma\ 
$$
does hold.}
\medskip
\noindent
{\bf Proof:} Replacing $V_{\sigma}$ with $V$ if
necessary, we may assume $\dim\sigma = n$. Take again 
a line $\ell$ intersecting $\buildrel \circ \over
\sigma$ and set $W:=V/\ell, B\b:=S\b (W^*)$ with
the quotient projection $\pi:V \longrightarrow W$, and, finally,
$\Lambda:=\pi (\partial \sigma)$. Choose
$T \in A^2$ with $T|_\ell \not= 0$. Then we have
$A\b \cong B\b[T]$ and
$$
\overline E\b_\sigma
\cong
\overline E\b_{\partial \sigma}
\cong
\coker (\mu_T: (B\b/B^{>0}) \otimes_{B\b} E\b_{\partial
\sigma} \longrightarrow
(B\b/B^{>0}) \otimes_{B\b} E\b_{\partial
\sigma})
$$
$$
\cong 
\coker
( L: IH\b (\Lambda) \longrightarrow
IH\b(\Lambda))\ ,
$$
where the Lefschetz operator $L$
is induced by multiplication with the strictly convex
piecewise linear function $\psi \in \A^2(\Lambda)$
being defined as
$$
\psi:=T \circ (\pi|_{\partial \sigma})^{-1}\ .
$$
By [Ka] we know that $L$ is onto in degrees
$q \ge n$, whence the first part of the statement 
follows. Since the pairing (DP) induces a dual pairing
$$
\overline
E\b_\Delta \times 
\overline
E\b_{(\Delta, \partial \Delta)}
\longrightarrow \bR\b[-2n]\ 
$$
of graded vector spaces, the second part is dual 
to the first one. Finally for a finitely generated
graded $A\b$-module $M\b$ one has $M^q=0$
for $q \le r$ if and only if
$\overline M^q=0$ for $q \le r$.
\qed
\bigskip
\noindent
{\bf Proof of 3.1}:
We have to show that for any cone $\sigma \not=o$
a homomorphism
$\phi_{\partial \sigma} \: E\b_{\partial \sigma}
\to F\b_{\partial \sigma}$, homogeneous of degree 0, 
extends in a unique
way to a homomorphism $\phi_{\sigma} \:
E\b_{\sigma} \to F\b_{\sigma}$. 
Since, according to 3.2 we have $E_{(\sigma, \partial \sigma)}^q 
=0=F_{(\sigma, \partial \sigma)}^q$ 
for $q \le \dim \sigma$,
the surjective restriction homomorphisms $E^q_\sigma \to
E^q_{\partial \sigma}$ and $F^q_\sigma \to
F^q_{\partial \sigma}$ are even isomorphisms in the same 
range. On the other hand, as a consequence of the 
first part of 3.2, the $A\b$-modules
$E\b_\sigma$ and $F\b_\sigma$ can be generated
by homogeneous elements of degree below
$\dim \sigma$, whence the assertion follows. \qed
\bigskip
As a consequence of 3.1 the duality isomorphism
$(DI)$ is unique up to a non-zero real factor.
And we
arrive at a natural choice 
when we fix an isomorphism
$E\b_o \cong \bR\b$.
In fact we require then 
$1 \in E\b_o$ to be mapped to the dual element
$1^* \in (E\b_o)^* = (D\E)\b_o$. So finally we may 
define the intersection product using the 
isomorphisms $(DI)$ and $D(\omega)$. 

Let us now compare with other possibilities
to define an intersection product. First we
review the simplicial case: There the 
intersection product
is nothing but the product of functions
$$
\A\b \times \A\b \longrightarrow \A\b
$$
taken at the level of global sections
$$
A\b_\Delta \times A\b_{(\Delta, \partial \Delta)} 
\longrightarrow 
A\b_{(\Delta, \partial \Delta)}
$$
followed by an evaluation map
$$
e_\Delta: A\b_{(\Delta, \partial \Delta)}
\longrightarrow A\b[-2n]\ .
$$
That strategy can also be applied in the general case,
since it is possible to define an (in fact non-canonical)
bilinear map of sheaves
$$
(.,.):\E\b \times \E\b \longrightarrow \E\b
$$
extending the multiplication
$$
E\b_o \times E\b_o
\longrightarrow E\b_o
$$
of real numbers.
Then again compose with an evaluation map
$$
e_\Delta:=e_\Delta (\omega): E\b_{(\Delta, \partial \Delta)}
\longrightarrow A\b[-2n]\ ,
$$ 
the image of
$1 \in E^0_o \cong  E^0_\Delta \subset 
E\b_\Delta$ in $(D\E)\b_\Delta
=\Hom (E\b_{(\Delta, \partial \Delta)},
A\b)
\otimes \omega$
\par
\noindent
$\cong 
\Hom (E\b_{(\Delta, \partial \Delta)},
 A\b[-2n])$
under the duality isomorphism $(DI)$.
Now, there is an associated sheaf homomorphism
$$
\E\b \longrightarrow D\E\ :
$$
For any cone $\sigma \in \Delta$ and $\omega_\sigma \in
\det (V_\sigma^*)$ we define
$e_\sigma (\omega_\sigma)$ in analogy to
$e_\Delta (\omega)$ and note that
$e_\sigma (\omega_\sigma) \otimes \omega_\sigma$
does not depend on the choice of $\omega_\sigma
\not=0$. Then, over $\sigma \in \Delta$, our sheaf 
homomorphism is 
$$
E\b_\sigma \longrightarrow (D\E)_\sigma,
f \mapsto (e_\sigma (\omega_\sigma)\circ (f,..) )
\otimes \omega_\sigma\ .
$$
In fact, over the zero cone it coincides with
our previous choice of an isomorphism 
$\E\b \longrightarrow D\E$
and thus both agree,
according to 3.2; so also the corresponding intersection products
coincide.

In [BBFK] such a bilinear sheaf homomorphism has recursively
been constructed on the $k$-skeletons of $\Delta$, an 
other possible choice is as follows:
Consider a simplicial refinement $\hat \Delta$ of $\Delta$
and denote $\hat \A\b$ the sheaf of piecewise polynomial functions on
$\hat \Delta$. Then, according to 1.5 we have
$$
\pi_*( \hat \A\b) \cong \E\b \oplus \K\b
$$
with some perverse sheaf $\K\b$. Here $\pi:=\id_V$.
Now take the following bilinear sheaf homomorphism
$$
\E\b \times \E\b \subset
\pi_*( \hat \A\b) \times \pi_*( \hat \A\b)
\buildrel \pi_*(\mu) \over
\longrightarrow
\pi_*( \hat \A\b) 
\buildrel p \over
\longrightarrow
\E\b\ ,
$$
where $\mu: \hat \A\b \times \hat \A\b \longrightarrow
\hat \A\b$ is the multiplication of functions and
$p: \E\b \oplus \K\b \longrightarrow \E\b$
the projection on the first factor.
Finally let us remark that 
writing
$$
\hat A\b_{(\Delta,
\partial \Delta)}=E\b_{(\Delta,
\partial \Delta)} \oplus
K\b_{(\Delta,
\partial \Delta)}
$$
one has
$$
e_\Delta (\omega) =
e_{\hat \Delta}(\omega)|_{E\b_{(\Delta,
\partial \Delta)}}\ ,
$$
such that we may also define the intersection 
product via
$$
E\b_\Delta \times E\b_{(\Delta,
\partial \Delta)}
\subset
\hat A\b_\Delta \times \hat A\b_{(\Delta,
\partial \Delta)}
\buildrel \mu \over
\longrightarrow
\hat A\b_{(\Delta,
\partial \Delta)} 
\buildrel e_{\hat \Delta}(\omega)
\over \longrightarrow
A\b[-2n]\ .
$$
\bigskip
\centerline{\bf References}
\medskip\noindent

\item{}[BBFK] {\sc G.~Barthel, J.-P.~Brasselet, 
K.-H.~Fieseler, L.~Kaup:} {\it Combinatorial Intersection Cohomolgy for Fans
}, T\^ohoku Math. J. {\bf 54} (2002), 1--41.

\item{}[BreLu$_1$] {\sc P.~Bressler, V.~Lunts:} {\it Intersection
cohomology on non-rational polytopes}, (pr)e-print 
{\tt math.AG/0002006}, 2000, to appear in Compositio Math.

\item{}[BreLu$_2$] {\sc P.~Bressler, V.~Lunts:} 
{\it Hard Lefschetz Theorem and
Hodge-Riemann Relations for Intersection Cohomology of
Nonrational Polytopes}, 
\- (pr)e-print 
{\tt math.AG/0302236v2},
2002.

\item{}[Bri] {\sc M.~Brion:} {\it The Structure of the Polytope
Algebra}, T\^ohoku Math. J. {\bf 49} (1997), 1--32.

\item{}[Ka] {\sc K.~Karu:} {\it Hard Lefschetz 
Theorem for Nonrational 
Polytopes}, (pr)e-print {\tt math.AG/0112087},
2003.
\bigskip
\baselineskip 2 mm
\item{} K.-H. Fieseler
\item{} Matematiska Institutionen
\item{} Box 480
\item{} Uppsala Universitet
\item{} SE-75106 Uppsala
\item{} e-mail: {\tt khf@math.uu.se}


\bye